\documentclass[a4paper,11pt]{article}
\usepackage{amsmath, amssymb, amsthm}
\usepackage{graphicx} 
\usepackage[all]{xy}
\usepackage[latin1]{inputenc}

\begin{document}

\def\Homeo{{\rm Homeo}} \def\MM{\mathcal M} \def\D{\mathbb D}
\def\C{\mathbb C} \def\S{\mathbb S} \def\Id{{\rm Id}}
\def\R{\mathbb R} \def\Z{\mathbb Z} \def\Ker{{\rm Ker}}
\def\CC{\mathcal C} \def\AA{\mathcal A} \def\N{\mathbb N}
\def\SSS{\mathfrak S} \def\ZZ{\mathcal Z}

%%%%%%%%%%%%%%%%%%%%%%%%%%%%%%%%%%%%%%%%%%%%%%%%%%%%%%%%%%%%%

\title{\bf{Mapping class groups of non-orientable surfaces for beginners}}

\author{\textsc{Luis Paris}}

\date{\today}

\maketitle

\begin{abstract}
\noindent
The present paper are the notes of a mini-course addressed mainly to non-experts.
It purpose it to provide a first approach to the theory of mapping class groups of non-orientable surfaces. 
\end{abstract}


\bigskip\bigskip\noindent
The present paper are the notes of a mini-course given in a workshop entitled ``Winter Braids'', held in Dijon, France, from February 10th to February 13th, 2014. 
It is addressed mainly to non-experts.
It purpose it to provide a first approach to the theory of mapping class groups, in general, and mapping class groups of non-orientable surfaces, in particular. 
We refer to \cite{FarMar1} for a detailed book on mapping class groups of orientable surfaces, and to \cite{Korkm1} for a survey on mapping class groups of non-orientable surfaces. 
There is no book on mapping class groups of non-orientable surfaces. 

\bigskip\noindent
{\bf Definition.}
Let $M$ be a surface, which may be orientable or not, have boundary or not, and have more than one connected component. 
Let $B$ be a union of boundary components of $M$ (which may be empty), and let $\Homeo(M,B)$ denote the group of homeomorphisms $h:M \to M$ that pointwise fix $B$. 
We endow $\Homeo(M,B)$ with the so-called \emph{compact-open topology}. 
Recall that this topology is defined as follows. 
For a compact subset $K$ of $M$ and an open subset $U$ we set $\Delta (K,U) = \{ h \in \Homeo (M,B) \mid h(K) \subset U\}$. 
Then the set of all $\Delta(K,U)$'s forms a basis for the topology of $\Homeo(M,B)$.
The following result is left as an exercise to the experts in topology. 
The others may just admit it. 

\bigskip\noindent
{\bf Theorem 1.}
{\it The group $\Homeo(M,B)$ is arcwise locally connected.}

\bigskip\noindent
{\bf Definition.}
We say that two homeomorphisms $f,g \in \Homeo(M,B)$ are \emph{isotopic} if they are in the same connected component of $\Homeo(M,B)$.
Let $\Homeo_0(M,B)$ denote the connected component of $\Homeo(M,B)$ that contains the identity. 
The \emph{mapping class group} of the pair $(M,B)$ is defined to be the quotient group $\MM(M,B)= \Homeo(M,B) /\Homeo_0(M,B)$. 
In other words, $\MM(M,B)$ is the group of isotopy classes of elements of $\Homeo(M,B)$. 
This is a discrete group. 
The group $\MM(M,\emptyset)$ (with $B=\emptyset$) will be simply denoted by $\MM(M)$.

\bigskip\noindent
{\bf Remark.}
Suppose that $M$ is orientable and connected.
Let $\Homeo^+(M,B)$ denote the subgroup of $\Homeo(M,B)$ consisting of the homeomorphisms that preserve the orientation. 
We have $\Homeo_0(M,B) \subset \Homeo^+(M,B)$, and $\MM^+(M,B) = \Homeo^+(M,B)/\Homeo_0(M,B)$ is a subgroup of $\MM(M,B)$.
In 99.9\% of the cases, the group called \emph{mapping class group} is $\MM^+(M,B)$ and not $\MM(M,B)$.
The latter is usually called the \emph{extended mapping class group}.
Actually, I know only one exception to this rule, a paper by Ivanov \cite{Ivano1}.

\bigskip\noindent
{\bf Exercise.}
Suppose that $M$ is a connected and orientable surface. 
\begin{itemize}
\item[(1)]
What is the meaning of being orientable for a homeomorphism of $M$?
\item[(2)]
What is the index of $\MM^+(M,B)$ in $\MM(M,B)$?
\end{itemize}

\bigskip\noindent
{\bf Definition.}
A \emph{crosscap} is a graphical representation of the M\"obius ribbon, as shown in Figure~1, where, in the circle which bounds the disk with a cross, two antipodal points are identified.  

\begin{figure}[tbh]
\bigskip\centerline{
\setlength{\unitlength}{0.4cm}
\begin{picture}(3,3)
\put(0,0){\includegraphics[width=1.2cm]{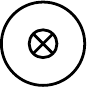}}
\end{picture}} 

\bigskip
\centerline{{\bf Figure 1.} A crosscap.}
\end{figure}

\bigskip\noindent
{\bf Graphical representation.}
Let $N$ be a non-orientable connected surface of genus $\rho$ with $m$ boundary components.
Pick $S$, an orientable connected surface of genus $\rho_1$ with $\rho_2+m$ boundary components, where $2\rho_1+\rho_2=\rho$, and $\rho_2 \ge 1$. 
Then we represent $N$ drawing crosses, that represent crosscaps, to $\rho_2$ of the boundary components of $S$ in a standard picture of $S$ (see Figure 2).

\bigskip\noindent
{\bf Example.}
In Figure 2 we have two different representations of the non-orientable surface of genus $3$ with $1$ boundary component. 

\begin{figure}[tbh]
\begin{center}
\begin{tabular}{ccc}
\parbox[c]{3.4cm}
{\setlength{\unitlength}{0.4cm}
\begin{picture}(8,4)
\put(0,0){\includegraphics[width=3cm]{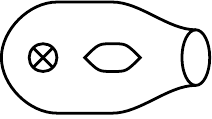}}
\end{picture}}
&
\parbox[c]{3.8cm}
{\setlength{\unitlength}{0.4cm}
\begin{picture}(9,4)
\put(0,0){\includegraphics[width=3.4cm]{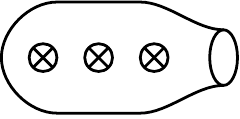}}
\end{picture}}
\end{tabular}

\bigskip
{\bf Figure 2.} Non-orientable surface of genus $3$ with $1$ boundary component. 
\end{center}
\end{figure}

\bigskip\noindent
The following is known as the ``Alexander theorem''. 

\bigskip\noindent
{\bf Proposition 2}
(Alexander theorem).
{\it Let $\D=\{ z \in \C; |z|\le 1\}$ be the standard disk. 
Then $\MM(\D,\partial \D) = \{1\}$.}

\bigskip\noindent
{\bf Proof.}
We need to prove that $\Homeo(\D,\partial \D)$ is connected. 
Let $h\in \Homeo (\D, \partial \D)$. 
For $t \in (0,1]$ we define $h_t \in \Homeo(\D,\partial \D)$ by:
\[
h_t(z)=\left\{\begin{array}{ll}
z &\quad\text{if } |z|\ge t\,,\\
t\,h(\frac{1}{t} z) &\quad\text{if } |z|<t\,,
\end{array}\right.
\]
and we set $h_0 = \Id$.
Then $\{h_t\}_{t \in [0,1]}$ is a path in $\Homeo(\D,\partial \D)$ joining the identity to $h$.
\qed

\bigskip\noindent
We turn now to give some useful definitions and results, often implicitly assumed in the literature. 

\bigskip\noindent
{\bf Definition.}
Let $X,Y$ be two topological spaces.
We say that two maps $f,g : X \to Y$ are \emph{homotopic} if there exists a continuous family $\{f_t: X \to Y\}_{t \in [0,1]}$  of continuous maps such that $f_0=f$ and $f_1=g$.
We say that two embeddings $f,g : X \to Y$ are \emph{isotopic} if there exists a continuous family $\{f_t : X \to Y\}_{t \in [0,1]}$ of embeddings such that $f_0=f$ and $f_1=g$.
Finally, we say that two embeddings $f,g : X \to Y$ are \emph{ambient isotopic} if there exists a continuous family $\{\psi_t : Y \to Y\}_{t \in [0,1]}$ of homeomorphisms such that $\psi_0=\Id$ and $\psi_1 \circ f = g$.
The families $\{f_t\}_{t \in [0,1]}$ and $\{ \psi_t\}_{t \in [0,1]}$ are both called \emph{isotopies}.
It is clear that ``ambient isotopic'' implies ``isotopic'' which implies ``homotopic''.
We will see that the converse holds in some cases. 

\bigskip\noindent
{\bf Definition.}
We come back to our surface $M$.
A \emph{circle} in $M$ is an embedding $a : \S^1 \to M$ of the standard circle whose image is in the interior of $M$. 
A circle $a$ is called \emph{essential} if it does not bound any disk embedded in $M$.
An \emph{arc} in $M$ is an embedding $a:[0,1] \to M$ of the standard interval such that $a(0), a(1) \in \partial M$ and $a((0,1)) \cap \partial M = \emptyset$.

\bigskip\noindent
{\bf Theorem 3}
(Epstein \cite{Epste1}).
{\it 
\begin{itemize}
\item[(1)]
Let $a_0,a_1 : \S^1 \to M$ be two essential circles.
Then $a_0$ and $a_1$ are homotopic if and only if they are ambient isotopic by an isotopy which pointwise fixes the boundary of $M$.
\item[(2)]
Let $a_0, a_1 : [0,1] \to M$ be two arcs having the same extremities.
Then $a_0$ and $a_1$ are homotopic with respect to the extremities if and only if $a_0$ and $a_1$ are ambient isotopic by an isotopy which pointwise fixes the boundary of $M$.
\end{itemize}}

\bigskip\noindent
Now, we turn to apply Theorem 3 to the proof of the following.

\bigskip\noindent
{\bf Proposition 4.}
{\it Let $C=\S^1 \times [0,1]$ be the cylinder, and let $c_0=\S^1\times \{0\}$ and $c_1 = \S^1 \times \{1\}$ be the boundary components of $C$ (see Figure 3\,(i)).
Then $\MM(C,c_0)=\{1\}$.}

\begin{figure}[tbh]
\begin{center}
\begin{tabular}{ccc}
\parbox[c]{2.2cm}
{\setlength{\unitlength}{0.4cm}
\begin{picture}(5,3)
\put(0.1,0){\includegraphics[width=1.6cm]{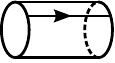}}
\put(0,2.3){\small $c_0$}
\put(2,0.8){\small $a$}
\put(3.5,2.3){\small $c_1$}
\end{picture}}
&
\parbox[c]{3.4cm}
{\setlength{\unitlength}{0.4cm}
\begin{picture}(8,5)
\put(0,1){\includegraphics[width=3.2cm]{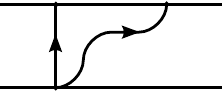}}
\put(1,2.5){\small $\tilde a$}
\put(1.7,0.4){\small $\tilde x_0$}
\put(4.5,2){\small $\tilde b$}
\put(5,4.6){\small $(n_0,1)$}
\end{picture}}
&
\parbox[c]{3.cm}
{\setlength{\unitlength}{0.4cm}
\begin{picture}(7,5)
\put(0,1){\includegraphics[width=2.8cm]{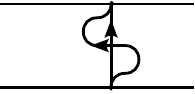}}
\put(2.4,2.8){\small $\tilde b$}
\put(3.5,0.4){\small $\tilde x_0$}
\put(4.4,3.2){\small $\tilde a$}
\end{picture}}
\\
(i) & (ii) & (iii)
\end{tabular}

\bigskip
{\bf Figure 3.} cylinder and arcs.
\end{center}
\end{figure}

\bigskip\noindent
{\bf Proof.}
We need to prove that $\Homeo(C,c_0)$ is connected.
Let $h \in \Homeo(C,c_0)$.
The restriction of $h^{-1}$ to $c_1$ is a homeomorphism $f' : c_1 \to c_1$ which preserves the orientation. 
Such a homeomorphism is isotopic to the identity, that is, there exists a continuous family $\{ f_r' : c_1 \to c_1 \}_{r \in [0,1]}$ of homeomorphisms such that $f_0'=\Id_{c_1}$ and $f_1'=f'$.
This isotopy can be extended to an isotopy $\{f_r: C \to C \}_{r \in [0,1]}$ such that $f_0=\Id_C$, and $f_r$ is the identity on $c_0$ for all $r$.
So, upon replacing $h$ by $f_1 \circ h$, we can assume that $h$ is the identity on $c_0 \sqcup c_1$.

\bigskip\noindent
Let $a : [0,1] \to C$ be the arc defined by $a(t) = (1,t)$ (see Figure 3\,(i)), and let $b=h \circ a$.
Let $p: \tilde C \to C$ be the universal cover of $C$.
Recall that $\tilde C = \R \times [0,1]$, and that $p$ is defined by $p(s,t)=(e^{2i\pi s},t)$.
We take $\tilde x_0=(0,0)$ as a basepoint for $\tilde C$, and we denote by $\tilde a$ and $\tilde b$ the lifts of $a$ and $b$ in $\tilde C$ of source $\tilde x_0$, respectively (see Figure 3\,(ii)).
We have $\tilde a(t) = (0,t)$ for all $t \in [0,1]$, and $\tilde b(1) \in \Z \times \{1\}$ (say $\tilde b(1) = (n_0,1)$).

\bigskip\noindent
For $r \in [0,1]$, we define $f_r : C \to C$ by 
\[
f_r(z,t) = (e^{-2i\pi rtn_0}\,z,t)\,.
\]
Then $\{f_r\}_{r \in [0,1]}$ is an isotopy which pointwise fixes $c_0$, and such that $f_0 = \Id_C$.
Let $b'=f_1 \circ h \circ a = f_1 \circ b$, and let $\tilde b'$ be the lift of $b'$ in $\tilde C$ of source $\tilde x_0$.
Then $\tilde b'(1) = \tilde a(1) = (0,1)$, and $f_1 \circ h$ is again the identity on $c_0 \sqcup c_1$.

\bigskip\noindent
So, upon replacing $h$ by $f_1 \circ h$, we can also assume that $\tilde b(1)=\tilde a(1)$ (see Figure~3\,(iii)).
Since the lifts of $a$ and $b$ of source $\tilde x_0$ have the same endpoint, the arcs $a$ and $b$ are homotopic relatively to they extremities. 
By Theorem 3 it follows that $a$ and $b$ are ambient isotopic by an isotopy which pointwise fixes the boundary of $C$.
In other words, there exists an isotopy $\{f_r: C \to C\}_{r \in [0,1]}$ such that $f_0=\Id_C$, $f_r$ is the identity on $c_0 \sqcup c_1$, and $f_1 \circ b = a$.

\bigskip\noindent
So, upon replacing $h$ by $f_1 \circ h$, we can assume that $h$ is the identity on $c_0 \cup c_1 \cup a$.
If we cut $C$ along $a$, we obtain a disk, and therefore we can apply the Alexander theorem (Proposition~2) to get an isotopy $\{f_r: C \to C\}_{r \in [0,1]}$ such that $f_0=\Id_C$, $f_1=h$, and $f_r$ is the identity on $c_0 \cup c_1 \cup a$.
\qed

\bigskip\noindent
{\bf Proposition 5}
(Epstein \cite{Epste1}).
{\it Let $N$ be the M\"obius ribbon.
Then $\MM(N, \partial N)=\{1\}$.}

\bigskip\noindent
The proof of Proposition 5 is based on Lemma 6 below. 
This is a surprise in the sense that there is no such a result of ``rigidity'' in the world of orientable surfaces.

\bigskip\noindent
We picture the M\"obius ribbon as in Figure 4\,(i), where $c_0$ and $c_1$ are identified according to the orientation shown in the figure.
We denote by $a$ the arc depicted in the figure, by $A$ its source, and by $B$ its target.
Note that $a$ does not separate $N$ into two connected components.
An example of an arc from $A$ to $B$ which separates $N$ into two connected components is depicted in Figure 4\,(ii).

\begin{figure}[tbh]
\begin{center}
\begin{tabular}{cc}
\parbox[c]{3.4cm}
{\setlength{\unitlength}{0.4cm}
\begin{picture}(8,5)
\put(0.8,1){\includegraphics[width=2.6cm]{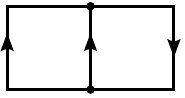}}
\put(0,2.5){\small $c_0$}
\put(3.2,2.5){\small $a$}
\put(3.5,0.3){\small $A$}
\put(3.5,4.5){\small $B$}
\put(7.3,2.5){\small $c_1$}
\end{picture}}
&
\parbox[c]{2.8cm}
{\setlength{\unitlength}{0.4cm}
\begin{picture}(6,5)
\put(0,1){\includegraphics[width=2.6cm]{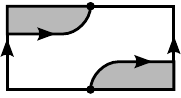}}
\put(2.8,0.3){\small $A$}
\put(2.8,4.5){\small $B$}
\end{picture}}
\\
(i) & (ii)
\end{tabular}

\bigskip
{\bf Figure 4.} M\"obius ribbon and arcs.
\end{center}
\end{figure}

\bigskip\noindent
{\bf Lemma 6}
(Epstein \cite{Epste1}).
{\it Let $b$ be an arc in the M\"obius ribbon $N$ from $A$ to $B$ which does not separate $N$ into two connected components. 
Then $b$ is homotopic to $a$ relatively to the extremities.}

\bigskip\noindent
{\bf Remark.}
By Theorem 3, the arcs $a$ and $b$ of Lemma 6 are therefore ambient isotopic by an isotopy which pointwise fixes the boundary of $N$. 

\bigskip\noindent
{\bf Proof.}
Let $\tilde N \to N$ be the universal cover of $N$, and let $\gamma : \tilde N \to \tilde N$ be the generator of its deck transformation group (see Figure 5\,(i)).
We choose a lift $\tilde A$ of $A$, we denote by $\tilde a$ the lift of $a$ of source $\tilde A$, and we denote by $\tilde B$ the endpoint of $\tilde a$.
Let $\tilde b$ be the lift of $b$ of source $\tilde A$.
There exists $m \in \Z$ such that the endpoint of $\tilde b$ is $\gamma^m \tilde B$.
If $m$ were odd with $m \ge 3$, then $\tilde b$ and $\gamma^2 \tilde b$ would intersect (see Figure 5\,(ii)), which is not possible.
Similarly, we cannot have $m$ odd with $m \le -3$.
If $m=1$ or $-1$, then $b$ would separate $N$ into two components (see Figure~4\,(ii)).
If $m$ were even with $m \ge 2$, then $\tilde b$ and $\gamma \tilde b$ would intersect (see Figure 5\,(iii)), which is not possible. 
Similarly, we cannot have $m$ even and $m \le -2$.
In conclusion, we have $m=0$, hence $b$ is homotopic to $a$ relatively to the extremities.
\qed

\begin{figure}[tbh]
\bigskip\centerline{
\setlength{\unitlength}{0.4cm}
\begin{picture}(24,4)
\put(0,1){\includegraphics[width=9.6cm]{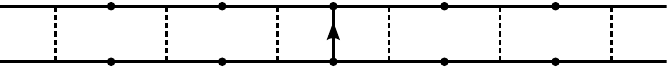}}
\put(3,0.1){\small $\gamma^2\,\tilde A$}
\put(3,3.5){\small $\gamma^2\,\tilde B$}
\put(7.3,0.1){\small $\gamma\, \tilde B$}
\put(7.3,3.5){\small $\gamma\,\tilde A$}
\put(11.5,0.1){\small $\tilde A$}
\put(11.5,3.5){\small $\tilde B$}
\put(15,0.1){\small $\gamma^{-1}\,\tilde B$}
\put(15,3.5){\small $\gamma^{-1}\,\tilde A$}
\put(19,0.1){\small $\gamma^{-2}\,\tilde A$}
\put(19,3.5){\small $\gamma^{-2}\,\tilde B$}
\put(11.2,2.3){\small $\tilde a$}
\end{picture}} 

\centerline{(i)}

\bigskip\bigskip\centerline{
\setlength{\unitlength}{0.4cm}
\begin{picture}(28,3)
\put(0,1){\includegraphics[width=11.2cm]{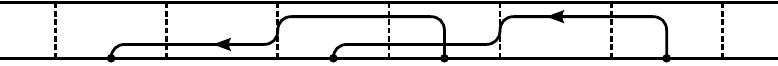}}
\put(6.5,2.1){\small $\gamma^2\,\tilde b$}
\put(11,0.1){\small $\gamma^3\,\tilde B$}
\put(20,1.5){\small $\tilde b$}
\put(23.5,0.1){\small $\tilde A$}
\end{picture}} 

\centerline{(ii)}

\bigskip\bigskip\centerline{
\setlength{\unitlength}{0.4cm}
\begin{picture}(20,4)
\put(0,1){\includegraphics[width=8cm]{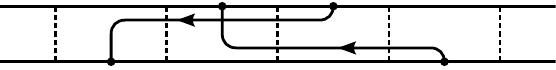}}
\put(6.3,1.7){\small $\gamma\,\tilde b$}
\put(7,3.6){\small $\gamma^2\,\tilde B$}
\put(13,2){\small $\tilde b$}
\put(15.5,0.1){\small $\tilde A$}
\end{picture}} 

\centerline{(iii)}

\bigskip
\centerline{{\bf Figure 5.} Universal cover of the M\"obius ribbon.}
\end{figure}

\bigskip\noindent
{\bf Proof of Proposition 5.}
Let $N$ be the M\"obius ribbon.
As ever, we need to show that $\Homeo(N,\partial N)$ is connected.
Let $h \in \Homeo(N, \partial N)$.
Consider the arc $a$ of Lemma 6, and set $b=h(a)$.
Then $b$ is an arc joining $A$ to $B$ which does not separate $N$ into two components, hence, by Lemma 6, $a$ and $b$ are ambient isotopic by an isotopy which pointwise fixes the boundary of $N$. 
In other words, there exists an isotopy $\{f_r: N \to N\}_{r \in [0,1]}$ such that $f_0=\Id_N$, $f_r$ is the identity on $\partial N$, and $f_1 \circ b = a$.

\bigskip\noindent
Upon replacing $h$ by $f_1 \circ h$, we may assume that $h$ is the identity on $\partial N \cup a$.
If we cut $N$ along $a$, we obtain a disk, and therefore we can apply the Alexander theorem (Proposition 2) to get an isotopy $\{f_r: N \to N\}_{r \in [0,1]}$ such that $f_0=\Id_N$, $f_1=h$, and $f_r$ is the identity on $\partial N \cup a$.
\qed

\bigskip\noindent
{\bf Definition.}
Recall that the cylinder is the surface $C=\S^1 \times [0,1]$.
The \emph{standard Dehn twist} is the homeomorphism $T: C \to C$ defined by
\[
T(z,t)= (e^{2i\pi t} z,t)\,.
\]
We also call \emph{standard Dehn twist} the mapping class $t \in \MM(C, \partial C)$ of $T$ (see Figure 6).

\begin{figure}[tbh]
\bigskip\centerline{
\setlength{\unitlength}{0.4cm}
\begin{picture}(12,3)
\put(0.5,0){\includegraphics[width=4.4cm]{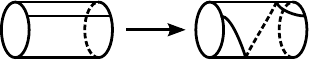}}
\put(5.5,1.5){\small $T$}
\end{picture}} 

\bigskip
\centerline{{\bf Figure 6.} Standard Dehn twist.}
\end{figure}

\bigskip\noindent
The main ingredients of the proof of the next proposition lie more or less in the proof of Proposition~4. 
The proof is left to the reader.

\bigskip\noindent
{\bf Proposition 7.}
{\it Let $C$ be the cylinder, and let $t \in \MM(C, \partial C)$ be the standard Dehn twist.
Then $\MM(C,\partial C)$ is an infinite cyclic group generated by $t$.}

\bigskip\noindent
{\bf Definition.}
Let $M$ be a surface, and let $B$ be a union of boundary components of $M$.
Let $M'$ be a subsurface of $M$, and let $B'$ be a union of boundary components of $M'$.
We say that the pair $(M',B')$ is \emph{admissible} if it satisfies the following condition.
\begin{itemize}
\item
If $a$ is a boundary component of $M'$, then either $a$ is a boundary component of $M$, or $\partial M \cap a = \emptyset$.
Moreover, if $a \cap \partial M = \emptyset$ or if $a \subset B$, then $a \subset B'$.
\end{itemize}
Let $(M',B')$ be an admissible pair.
Then there is a natural injective homomorphism $\iota : \Homeo(M',B') \to \Homeo(M,B)$.
The image of an $F \in \Homeo(M',B')$ under this homomorphism is obtained by extending $F$ to $M\setminus M'$ by the identity.
This homomorphism induces a homomorphism $\iota_*: \MM(M',B') \to \MM(M,B)$.
An example of an admissible pair is given in Figure 7.

\begin{figure}[tbh]
\bigskip\centerline{
\setlength{\unitlength}{0.4cm}
\begin{picture}(11,8)
\put(0,1){\includegraphics[width=4.2cm]{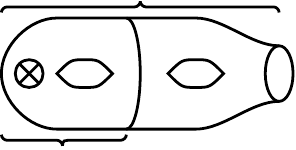}}
\put(1.7,0.3){\small $M'$}
\put(4.5,6.6){\small $M$}
\end{picture}} 

\bigskip
\centerline{{\bf Figure 7.} Admissible subsurface.}
\end{figure}

\bigskip\noindent
The above defined homomorphism $\iota_*: \MM(M',B') \to \MM(M,B)$ is not always injective, but its kernel is pretty well-understood.
An explicit description of $\Ker\, \iota_*$ may be found in \cite{ParRol1} for $M$ connected and orientable, and in \cite{Stuko1} for $M$ connected and non-orientable.
The general case is easily deduced from the connected case.
Here is an example of the type of result that we have. 

\bigskip\noindent
{\bf Theorem 8}
(Paris, Rolfsen \cite{ParRol1}, Stukow \cite{Stuko1}).
{\it Let $M$ be a connected surface, and let $(M',B')$ be an admissible pair.
Then $\iota_*: \MM(M',B') \to \MM(M)$ is injective if and only if the connected components of $M \setminus M'$ have all negative Euler characteristics.}

\bigskip\noindent
{\bf Definition.}
Recall that a \emph{circle} in $M$ is an embedding $a : \S^1 \to M$ of the standard circle whose image lies in the interior of $M$.
We assume the circles non-oriented, and we often identify them with their images.
Recall that a circle is called \emph{essential} if it does not bound any disk embedded in $M$.
We say that an essential circle $a$ is \emph{generic} if it does not bound any Möbius ribbon, and if it is not isotopic to any boundary component.
The isotopy class of a circle $a$ will be denoted by $[a]$.
Depending on whether a regular neighborhood of $a$ is an annulus or a Möbius ribbon, we will say that $a$ is \emph{two-sided} or \emph{one-sided}.
We denote by $\CC(M)$ the set of isotopy classes of generic circles, and by $\CC_0(M)$ the set of isotopy classes of two-sided generic circles. 

\bigskip\noindent
{\bf Definition.}
The \emph{intersection index} of two classes $\alpha, \beta \in \CC(M)$ is $i(\alpha, \beta) = \min\{|a \cap b| \mid a \in \alpha \text{ and } b \in \beta\}$.

\bigskip\noindent
{\bf Definition.}
Let $a: \S^1 \hookrightarrow \Sigma$ be a two-sided circle in $M$.
Take a tubular neighborhood of $a$, that is, an embedding $\varphi: C \hookrightarrow M$ of the annulus $C=\S^1 \times [0,1]$, such that:
\begin{itemize}
\item
$\varphi (C) \cap \partial M = \emptyset$,
\item
$\varphi(z,\frac{1}{2}) = a(z)$ for all $z \in \S^1$.
\end{itemize}
Let $T_a$ be the element of $\Homeo(M)$ defined by: 
\begin{itemize}
\item
$T_a(x) = (\varphi \circ T \circ \varphi^{-1}) (x)$ if $x$ is in the image of $\varphi$, where $T$ denotes the standard Dehn twist, 
\item
$T_a(x)=x$ if $x$ does not belong to the image of $\varphi$.
\end{itemize}
The isotopy class of $T_a$, denoted by $t_a \in \MM(M)$, is called the \emph{Dehn twist} along $a$ (see Figure 8).

\begin{figure}[tbh]
\bigskip\centerline{
\setlength{\unitlength}{0.4cm}
\begin{picture}(17,4)
\put(0,0){\includegraphics[width=6.6cm]{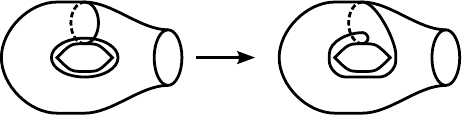}}
\put(2.5,0.6){\small $b$}
\put(3.6,3){\small $a$}
\put(7.5,2.5){\small $t_a$}
\put(11.5,0.5){\small $t_a(b)$}
\end{picture}} 

\bigskip
\centerline{{\bf Figure 8.} A Dehn twist.}
\end{figure}

\bigskip\noindent
{\bf Remark.}
Recall that, by Proposition 7, $\MM(C, \partial C)$ is an infinite cyclic group generated by the standard Dehn twist, $t$.
Recall also that the embedding $\varphi : C \hookrightarrow M$ induces a homomorphism $\varphi_* : \MM(C, \partial C) \to \MM(M)$.
Then $t_a = \varphi_*(t)$.

\bigskip\noindent
{\bf Proposition 9.}
{\it \begin{itemize}
\item[(1)]
The definition of $t_a$ depends neither on the choice of the tubular neighborhood, nor on the choice of an orientation of $a$.
However, it depends on the choice of an orientation on a regular neighborhood of $a$ (i.e. the image of $C$). 
\item[(2)]
If $a$ and $b$ are isotopic circles, then $t_a=t_b$.
\end{itemize}}

\bigskip\noindent
{\bf Definition.}
Let $\alpha \in \CC_0(M)$.
We choose $a \in \alpha$, and we define the \emph{Dehn twist} along $\alpha$, denoted by $t_\alpha$, as being the Dehn twist along $a$, that is, $t_\alpha = t_a$.
If $M$ is an orientable connected surface, then we may fix an orientation on $M$, and therefore the orientation of any regular neighborhood of $a$ is naturally well-defined, hence the element $t_a=t_\alpha$ is well-defined, as well. 
However, if $M$ is a non-orientable connected surface, then $t_a=t_\alpha$ depends on the choice of an orientation on a regular neighborhood of $a$.
We use to indicate in the figures this chosen local orientation, that is, the ``direction'' of the Dehn twist, with arrows like in Figure 9.

\begin{figure}[tbh]
\bigskip\centerline{
\setlength{\unitlength}{0.4cm}
\begin{picture}(20,6)
\put(0,0){\includegraphics[width=8cm]{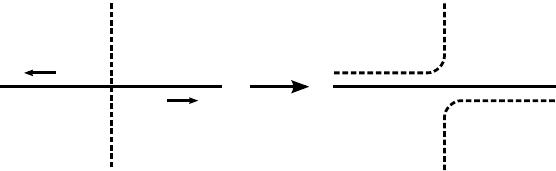}}
\put(6,3.2){\small $a$}
\put(9.5,3.3){\small $t_a$}
\put(18,3.2){\small $a$}
\end{picture}} 

\bigskip
\centerline{{\bf Figure 9.} Dehn twist on a non-orientable surface.}
\end{figure}

\bigskip\noindent
{\bf Remark.}
The notion of a Dehn twist along a one-sided circle does not exist. 

\bigskip\noindent
When $M$ is an orientable connected surface, the Dehn twists and the intersection indices are linked by the following formulas. 
These play an important role in the study of the mapping class groups of orientable surfaces. 

\bigskip\noindent
{\bf Theorem 10}
(Fathi, Laudenbach, Poénaru \cite{FaLaPo1}).
{\it
Assume that $M$ is an oriented connected surface. 
\begin{itemize}
\item[(1)]
Let $\alpha, \beta$ be two elements of $\CC(M)$, and let $k$ be an integer.
Then
\[
i(t_\alpha^k(\beta), \beta) = |k|\,i(\alpha, \beta)^2\,.
\]
\item[(2)]
Let $\alpha, \beta, \gamma$ be three elements of $\CC(M)$, and let $k$ be an integer.
Then
\[
\big|i(t_\alpha^k(\beta), \gamma)|- |k|\, i(\alpha, \beta)\, i(\alpha, \gamma) \big| \le i(\beta, \gamma)\,.
\]
\end{itemize}}

\bigskip\noindent
The above formulas are false if $M$ is a non-orientable surface (see \cite{Stuko2,Stuko3}).
There are other formulas for the non-orientable surfaces (see \cite{Stuko2}), that are less powerful than the ones for the orientable surfaces, but that are still useful.
Here are some. 

\bigskip\noindent
{\bf Theorem 11}
(Stukow \cite{Stuko2}).
{\it Assume that $M$ is a connected non-orientable surface.
Let $\alpha, \beta$ be two elements of $\CC_0(M)$, and let $k$ be a non-zero integer.
\begin{itemize}
\item[(1)]
$i(t_\alpha^k(\beta), \beta)=|k|$ if $i(\alpha, \beta)=1$.
\item[(2)]
$i(t_\alpha^k(\beta), \beta) \ge i(\alpha, \beta)$.
\item[(3)]
$i(t_\alpha^k (\beta), \beta) \ge (|k|-1)\,i(\alpha,\beta)^2 +2\,i(\alpha,\beta) -2$.
\end{itemize}
In particular, if $i(\alpha, \beta) \neq 0$, then $i(t_\alpha^k(\beta), \beta) >0$.}

\bigskip\noindent
{\bf Lemma 12.}
{\it Let $a$ be an essential circle in $M$.
If either $a$ bounds a M\"obius ribbon, or $a$ is isotopic to a boundary component, then $t_a=1$.}

\bigskip\noindent
{\bf Proof.}
Assume that $a$ bounds a M\"obius ribbon.
The case where $a$ is isotopic to a boundary component can be proved in the same way. 
We denote again by $C=\S^1 \times [0,1]$ the standard annulus, and we assume given a tubular neighborhood $\varphi : C \to M$ of $a$.
There exists a subsurface $N$ of $M$ homeomorphic to the M\"obius ribbon and which contains $\varphi (C)$ in its interior. 
We denote by $\varphi' : C \to N$ the tubular neighborhood of $a$ given by restriction of $\varphi$ to $N$, and by $\iota : N \to M$ the inclusion map. 
We have $\varphi = \iota \circ \varphi'$, hence $\varphi_* = \iota_* \circ \varphi'_*$.
Moreover, $\varphi'_*(t)=1$ since, by Lemma~6, $\MM(N, \partial N) = \{1\}$.
So, $t_a = \varphi_*(t)=1$.
\qed

\bigskip\noindent
However, if $a$ is a generic circle (that is, $a$ is essential, it does not bound any M\"obius ribbon, and it is not isotopic to any boundary component), then we have the following. 

\bigskip\noindent
{\bf Proposition 13.}
{\it Let $a$ be a generic circle in $M$.
Then the Dehn twist $t_a$ is of infinite order.}

\bigskip\noindent
{\bf Sketch of the proof.}
We set $\alpha = [a]$.
We leave to the reader to show that there exists $\beta \in \CC_0(M)$ such that $i(\alpha, \beta) \neq 0$.
Let $k$ be a non-zero integer. 
By Theorems 10 and 11, we have $i(t_\alpha^k(\beta), \beta) \neq 0$, hence $t_\alpha^k(\beta)\neq \beta$, and therefore $t_\alpha^k \neq \Id$.
\qed

\bigskip\noindent
We turn now to explain the relations that hold between two given Dehn twists. 
We start with an obvious case.

\bigskip\noindent
{\bf Lemma 14.}
{\it Let $\alpha, \beta \in \CC_0(M)$.
If $i(\alpha, \beta)=0$, then $t_\alpha t_\beta = t_\beta t_\alpha$.}

\bigskip\noindent
In order to handle the case of two classes with intersection index $1$, we need to choose local orientations as follows. 
Let $\alpha, \beta$ be elements of $\CC_0(M)$ such that $i(\alpha, \beta) =1$.
We choose $a \in \alpha$ and $b \in \beta$ such that $|a \cap b|=1$.
Let $P$ be a regular neighborhood of $a \cup b$.
Then $P$ is an orientable surface of genus $1$ with one boundary component (see Figure 10). 
In the next lemma we will assume given an orientation on $P$, and we will assume that the given orientations on the regular neighborhoods of $a$ and $b$ are those induced by the one of $P$. 

\begin{figure}[tbh]
\bigskip\centerline{
\setlength{\unitlength}{0.4cm}
\begin{picture}(7,4)
\put(0,0){\includegraphics[width=2.6cm]{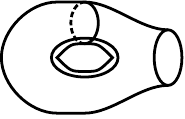}}
\put(2,0.7){\small $b$}
\put(3.6,3){\small $a$}
\end{picture}} 

\bigskip
\centerline{{\bf Figure 10.} Orientation of Dehn twists, case where $i(\alpha, \beta)=1$.}
\end{figure}

\bigskip\noindent
{\bf Lemma 15.}
{\it Let $\alpha, \beta \in \CC_0(M)$.
If $i(\alpha, \beta)=1$, then $t_\alpha t_\beta t_\alpha = t_\beta t_\alpha t_\beta$.}

\bigskip\noindent
{\bf Proof.}
Let $a \in \alpha$ and $b \in \beta$ such that $|a \cap b|=1$.  
We have
\[
t_\alpha t_\beta t_\alpha = t_\beta t_\alpha t_\beta \quad \Leftrightarrow \quad t_\beta^{-1} t_\alpha t_\beta = t_\alpha t_\beta t_\alpha^{-1} \quad \Leftrightarrow \quad t_{t_\beta^{-1}(\alpha)} = t_{t_\alpha(\beta)}\,.
\]
So, in order to show that $t_\alpha t_\beta t_\alpha = t_\beta t_\alpha t_\beta$, it suffices to show that $t_b^{-1}(a)$ is isotopic to $t_a(b)$.
As explained before, a regular neighborhood of $a \cup b$ is a one-holed torus, and it is easily seen that $t_a(b) \simeq t_b^{-1}(a)$ is the circle depicted in Figure 8, up to isotopy.
\qed

\bigskip\noindent
The proof of the following is more complicated than the proofs of Lemmas 14 and 15, especially for non-orientable surfaces.

\bigskip\noindent
{\bf Theorem 16}
(Ishida \cite{Ishid1}, Stukow \cite{Stuko3}).
{\it Let $\alpha, \beta \in \CC_0(M)$.
If $i(\alpha, \beta)\ge 2$, then $t_\alpha$ and $t_\beta$ freely generate a free group.}

\bigskip\noindent
{\bf Notation.}
From now on, we denote by $S_{\rho,m}$ the orientable surface of genus $\rho$ with $m$ boundary components, and by $N_{\rho,m}$ the non-orientable surface of genus $\rho$ with $m$ boundary components. 

\bigskip\noindent
{\bf Definition.}
Assume that $N=N_{1,2}$ is the one-holed M\"obius ribbon.
Let $a,c$ denote the boundary components of $N_{1,2}$, as shown in Figure 11. 
Let $s_c$ denote the mapping class in $\MM(N_{1,2},c)$ obtained by sliding the boundary component $a$ once along the path $u$ depicted in Figure 11.
Then $s_c$ is called the \emph{boundary slide} along $u$.

\begin{figure}[tbh]
\bigskip\centerline{
\setlength{\unitlength}{0.4cm}
\begin{picture}(7,5)
\put(0,0){\includegraphics[width=2.4cm]{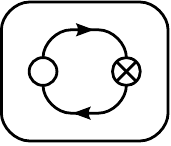}}
\put(0.7,3){\small $a$}
\put(3.5,4.2){\small $u$}
\put(6.1,3){\small $c$}
\end{picture}} 

\bigskip
\centerline{{\bf Figure 11.} Boundary slide.}
\end{figure}

\bigskip\noindent
{\bf Proposition 17}
(Korkmaz \cite{Korkm1}).
{\it The group $\MM(N_{1,2},c)$ is an infinite cyclic group generated by $s_c$.
The mapping class $s_c$ reverses the orientation in a regular neighborhood of $a$, and we have $s_c^2=t_c$, where $t_c$ denotes the Dehn twist along a circle in $N_{1,2}$ isotopic to $c$.}

\bigskip\noindent
{\bf Remark.}
By Lemma 12, we have $t_c=1$ in $\MM(N_{1,2})$.
However, $t_c$ is of infinite order in $\MM(N_{1,2},c)$.

\bigskip\noindent
More generally, we have the following definition. 

\bigskip\noindent
{\bf Definition.}
Let $a$ be a boundary component of $M$, and let $u$ be an oriented arc in $M$ whose extremities lie in $a$.
Let $s_{u,a}$ denote the mapping class in $\MM(M)$ obtained by sliding the boundary component $a$ once along the arc $u$.
Then $s_{u,a}$ is called the \emph{boundary slide}  of $a$ along $u$.

\bigskip\noindent
The proof of the following is left to the reader. 

\bigskip\noindent
{\bf Lemma 18.}
{\it If $u$ is two-sided, then $s_{u,a}$ is the product of two Dehn twists, $t_{b_1}$ and $t_{b_2}$, oriented according to Figure 12.}

\begin{figure}[tbh]
\bigskip\centerline{
\setlength{\unitlength}{0.4cm}
\begin{picture}(8,4)
\put(0.5,0){\includegraphics[width=2.8cm]{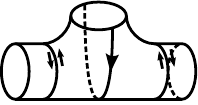}}
\put(1.5,2.5){\small $b_1$}
\put(4.5,0.5){\small $u$}
\put(2.7,3.5){\small $a$}
\put(5.8,2.5){\small $b_2$}
\end{picture}} 

\bigskip
\centerline{{\bf Figure 12.} Boundary slide along a two-sided arc.}
\end{figure}

\bigskip\noindent
By increasing the genus, that is, passing from genus $1$ to genus $2$, we meet new phenomena of ``rigidity'', as we turn now to show (see Propositions 19 and 21). 

\bigskip\noindent
{\bf Proposition 19}
(Lickorish \cite{Licko1}).
{\it We represent the Klein bottle, $N=N_{2,0}$, by two crosscaps on the sphere as depicted in Figure 13\,(i).
Let $\alpha$ be the class in $\CC_0(N)$ of the circle $a$ drawn in the figure.
Then $\CC_0(N)=\{\alpha\}$.
In other words, every two-sided generic circle in $N_{2,0}$ is isotopic to $a$.}

\begin{figure}[tbh]
\begin{center}
\begin{tabular}{cc}
\parbox[c]{1.8cm}
{\setlength{\unitlength}{0.4cm}
\begin{picture}(4,4)
\put(0,0){\includegraphics[width=1.6cm]{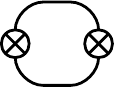}}
\put(1.9,3.3){\small $a$}
\end{picture}}
&
\parbox[c]{2.2cm}
{\setlength{\unitlength}{0.4cm}
\begin{picture}(5,4)
\put(0.5,0){\includegraphics[width=1.8cm]{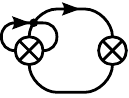}}
\put(1.2,3){\small $P$}
\put(2.5,1.4){\small $e$}
\put(3.2,3.2){\small $a$}
\end{picture}}
\\
\\
(i) & (ii)
\end{tabular}

\bigskip
{\bf Figure 13.} The Klein bottle.
\end{center}
\end{figure}

\bigskip\noindent
The proof of Proposition 19 given here is totally different from Lickorish's one \cite{Licko1}, and, as far as I know, is not written anywhere else in the literature.
Like the proof of Lemma 6, it is based on the study of the lifts of the (generic two-sided) circles to the universal cover of $N_{2,0}$.

\bigskip\noindent
{\bf Proof.}
We pick a generic two-sided circle $b$ of $N$, and we turn to prove that $b$ is isotopic to $a$.

\bigskip\noindent
{\bf Claim 1.}
{\it There exists a homeomorphism $H: N \to N$ that sends $b$ to $a$.}

\bigskip\noindent
{\bf Proof of Claim 1.}
Let $N'$ be the surface obtained by cutting $N$ along $b$.
Suppose first that $N'$ has two connected components, $N_1'$ and $N_2'$.
Since $b$ is generic and $N_i'$ is a subsurface of $N$, $N_i'$ has a unique boundary component, the genus of $N_i'$ is (strictly) less than $1$ if $N_i'$ is orientable, and the genus of $N_i'$ is (strictly) less than $2$ if $N_i'$ is non-orientable.
In other words, $N_i'$ is either a disk or a M\"obius ribbon: a contradiction. 
So, $N'$ is connected, $N'$ has two boundary components, and $\chi(N') = \chi(N)=0$.
This implies that $N'$ is an annulus.
In conclusion, $b$ has the same topological type as $a$, hence there exists a homeomorphism $H: N \to N$ that sends $b$ to $a$.
\qed

\bigskip\noindent
Let $c$ be a circle in any connected surface $M$.
We parametrize $c$, $c: [0,1] \to M$, we choose a base-point $P$ in $M$, we choose a path $d : [0,1] \to M$ connecting $P$ to $c(0)=c(1)$, and we denote by $\gamma$ the element of $\pi_1(M,P)$ represented by $d^{-1} c d$.
We say that $c$ is \emph{primitive} if $\gamma$ is not a non-trivial power in $\pi_1(M,P)$, that is, if $\mu \in \pi_1(M,P)$ and $k \in \Z$ are such that $\gamma= \mu^k$, then $k = \pm 1$.
It is easily seen that this definition depends neither on the choice of the parametrization of $c$, nor on the choice of $P$, nor on the choice of $d$.

\bigskip\noindent
{\bf Claim 2.}
{\it $b$ is primitive.}

\bigskip\noindent
{\bf Proof of Claim 2.}
By Claim 1, it suffices to prove Claim 2 for $b=a$.
Consider the base-point $P$, the orientation of $a$, and the loop $e$, depicted in Figure~13\,(ii).
We denote by $\alpha$ and $\epsilon$ the elements of $\pi_1(N,P)$ represented by $a$ and $e$, respectively. 
Then $\pi_1(N,P)$ has the following presentation.
\[
\pi_1(N,P)=\langle \alpha, \epsilon \mid \epsilon \alpha \epsilon^{-1} = \alpha ^{-1} \rangle \,.
\]
From this presentation follows that there exists a homomorphism $\varphi: \pi_1(N,P) \to \Z$ which sends $\alpha$ to $0$, and $\epsilon$ to $1$.
The kernel of $\varphi$ is an infinite cyclic group generated by $\alpha$.
If $\alpha$ were a power of some element $\beta$ of $\pi_1(N,P)$, then $\beta$ would belong to $\Ker\,\varphi = \langle \alpha \rangle$, hence $\alpha = \beta^{\pm 1}$.
\qed

\bigskip\noindent
{\bf Exercise.}
Prove that any generic circle in any connected surface is primitive. 

\bigskip\noindent
Consider the transformations $\sigma, \tau : \C \to \C$ defined as follows.
\[
\sigma(x+iy) = (x+1) +iy\,,\quad \tau(x+iy) = -x +i(y+1)\,.
\]
We can and do assume that $N$ is the quotient $N=\C/\langle \sigma, \tau \rangle$, and we denote by $p : \C \to N$ the quotient map.
Note that $p$ is the universal cover of $N$. 
We assume, furthermore, that $a$ is the circle parametrized by $a: [0,1] \to N$, $t \mapsto p(t) = p(t+0\,i)$. 

\bigskip\noindent
We compactify $\C$ with the asymptotic directions.
In other words, we set $\bar \C = \C \sqcup \partial \C$, where $\partial \C = \S^1 = \{\hat z \in \C \mid |\hat z|=1 \}$.
A sequence $\{z_n\}_{n \in \N}$ in $\C$ converges to $\hat z \in \partial \C$ if $\lim_{n \to \infty} |z_n| = \infty$ and $\lim_{n \to \infty} \frac{z_n}{|z_n|} = \hat z$.
Note that $\bar \C$ is homeomorphic to the disk $\D=\{z \in \C \mid |z| \le 1\}$ by a homeomorphism which sends $z$ to $\frac{z}{|z|+1}$ if $z \in \C$, and sends $\hat z$ to itself if $\hat z \in \partial \C = \S^1$.

\bigskip\noindent
We choose a parametrization of $b$, $b:[0,1] \to N$.
We can and do assume that $b(0) = b(1)= p(0)$.
We denote by $\tilde b : [0,1] \to \C$ the lift of $b$ in $\C$ of source $0$.
Then there exist $n,m \in \Z$ such that $\tilde b(1)=n +im$.
Note that $m$ must be even, since $b$ is two-sided.
Note also that either $n \neq 0$, or $m \neq 0$, since $b$ is not homotopic to the constant loop. 

\bigskip\noindent
{\bf Claim 3.}
{\it Either $n=0$ or $m=0$.}

\bigskip\noindent
{\bf Proof of Claim 3.}
Suppose that $n \neq 0$ and $m \neq 0$.
Set $\phi = \sigma^n\tau^m $.
Set $\tilde B = \cup_{ k \in \Z} \phi^k\, \tilde b$.
Observe that $\tilde B$ is a connected component of $p^{-1}(b)$ and is an arc in $\bar \C$ connecting $u_n = \frac{n+im}{\sqrt{n^2+m^2}} \in \partial \C$ to $-u_n \in \partial \C$ (see Figure 14).
On the other hand, $\tau\tilde B$ is another connected component of $p^{-1}(b)$, and $\tau\tilde B$ is an arc in $\bar \C$ connecting $v_n = \frac{-n+im}{\sqrt{n^2+m^2}} \in \partial \C$ to $-v_n \in \partial \C$.
But then $\tau\tilde B$ intersects $\tilde B$: a contradiction.
So, either $n=0$ or $m=0$.
\qed

\begin{figure}[tbh]
\bigskip\centerline{
\setlength{\unitlength}{0.4cm}
\begin{picture}(18,13)
\put(0,0){\includegraphics[width=7.2cm]{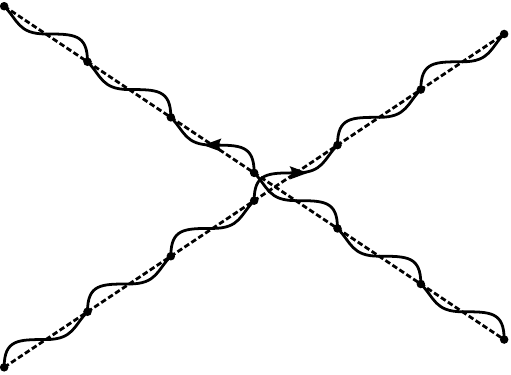}}
\put(0,11){\small $\tau \tilde B$}
\put(0.2,8.5){\small $-n+i(m+1)$}
\put(7,8.5){\small $\tau \tilde b$}
\put(10,7.5){\small $\tilde b$}
\put(12.3,7.7){\small $n+im$}
\put(17,10){\small $\tilde B$}
\end{picture}} 

\bigskip
\centerline{{\bf Figure 14.} Universal cover of $N_{2,0}$.}
\end{figure}

\bigskip\noindent
{\bf End of the proof of Proposition 19.}
By Claim 3, either $n=0$ or $m=0$.
If $n=0$, then, $m = \pm 1$ since by Claim 2 $b$ is primitive. 
But, as pointed out before, $m$ is even since $b$ is two-sided: a contradiction. 
Hence, $m=0$ and, since $b$ is primitive, $n=\pm 1$.
It follows that $b$ is homotopic to $a$, and therefore, by Theorem 3, $b$ is isotopic to $a$.
\qed

\bigskip\noindent
{\bf Exercise.}
Prove that there are precisely two isotopy classes of generic one-sided circles in the Klein bottle. 

\bigskip\noindent
Proposition 19 is the main ingredient in the proof of the following. 

\bigskip\noindent
{\bf Proposition 20}
(Lickorish \cite{Licko1}).
{\it The mapping class group of the Klein bottle, $N_{2,0}$, is isomorphic to $\Z/2\Z \times \Z /2\Z$.}

\bigskip\noindent
Proposition 19 extends to the one-holed Klein bottle as follows. 

\bigskip\noindent
{\bf Proposition 21}
(Stukow \cite{Stuko2}).
{\it We represent the one-holed Klein bottle, $N=N_{2,1}$, by two crosscaps on the disk as depicted in Figure 15\,(i).
Let $\alpha$ be the class in $\CC_0(N)$ of the circle $a$ drawn in the figure.
Then $\CC_0(N)=\{\alpha\}$.
In other words, every two-sided generic circle in $N_{2,1}$ is isotopic to $a$.}

\begin{figure}[tbh]
\begin{center}
\begin{tabular}{cc}
\parbox[c]{2.6cm}
{\setlength{\unitlength}{0.4cm}
\begin{picture}(6,5)
\put(0,0){\includegraphics[width=2.4cm]{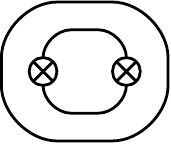}}
\put(2.8,4.2){\small $a$}
\end{picture}}
&
\parbox[c]{7.4cm}
{\setlength{\unitlength}{0.4cm}
\begin{picture}(18,3)
\put(0,0){\includegraphics[width=7.2cm]{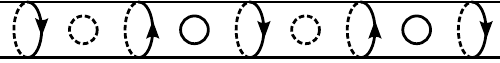}}
\put(0.5,2.5){\small $\tilde a_{-1}$}
\put(3.8,0.8){\small $\tilde d_{-1}$}
\put(4.5,2.5){\small $\tilde a_0$}
\put(7.5,0.5){\small $\tilde d_0$}
\put(8.5,2.5){\small $\tilde a_1$}
\put(11.5,0.5){\small $\tilde d_1$}
\put(12.5,2.5){\small $\tilde a_2$}
\put(15.5,0.5){\small $\tilde d_2$}
\put(16.5,2.5){\small $\tilde a_3$}
\end{picture}}
\\
(i) & (ii)
\end{tabular}

\bigskip
\begin{tabular}{ccc}
\parbox[c]{3cm}
{\setlength{\unitlength}{0.4cm}
\begin{picture}(7,4)
\put(0,1){\includegraphics[width=2.4cm]{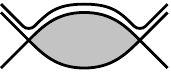}}
\put(0.5,3.2){\small $b'$}
\put(2.5,1.8){\small $\Delta$}
\put(6,0.6){\small $a$}
\put(6,2.4){\small $b$}
\end{picture}}
&
\parbox[c]{3cm}
{\setlength{\unitlength}{0.4cm}
\begin{picture}(7,5)
\put(1,1){\includegraphics[width=2cm]{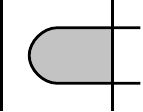}}
\put(0.8,0){\small $\tilde a_{k_0}$}
\put(3,2.7){\small $\Delta_0$}
\put(4.5,0){\small $\tilde a_{k_0+1}$}
\put(6.2,3.7){\small $\tilde b_0$}
\end{picture}}
&
\parbox[c]{3cm}
{\setlength{\unitlength}{0.4cm}
\begin{picture}(7,5)
\put(1,1){\includegraphics[width=2cm]{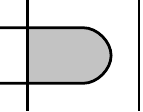}}
\put(0,3.7){\small $\tilde b_0$}
\put(1.5,0){\small $\tilde a_{k_1}$}
\put(2.9,2.7){\small $\Delta_1$}
\put(5.5,0){\small $\tilde a_{k_1+1}$}
\end{picture}}
\\
\\
(iii) & (iv) & (v)
\end{tabular}

\bigskip
{\bf Figure 15.} One-holed Klein bottle and bigons. 
\end{center}
\end{figure}

\bigskip\noindent
The proof of Proposition 21 that we give here is new and completely different from the one given in \cite{Stuko2}.
It is in the spirit of the proofs of Lemma 6 and Proposition 19.

\bigskip\noindent
{\bf Proof.}
Let $d$ denote the boundary component of $N_{2,1}$.
We assume that $N_{2,0}$ is obtained from $N_{2,1}$ by gluing a disk along $d$.
We denote by $\iota : N_{2,1} \to N_{2,0}$ the inclusion map, and by $\iota_* : \pi_1(N_{2,1}) \to \pi_1(N_{2,0})$ the homomorphism induced by $\iota$.
Note that $\iota_*$ is onto.

\bigskip\noindent
Recall the homomorphism $\varphi : \pi_1(N_{2,0}) \to \Z$ defined in the proof of Proposition 19.
Recall that this homomorphism sends $\alpha$ to $0$, and $\epsilon$ to $1$.
We set $\tilde \varphi = \varphi \circ \iota_* : \pi_1(N_{2,1}) \to \Z$, and we denote by $p : \tilde N \to N_{2,1}$ the associated regular cover.
We can describe $\tilde N$ as follows.
Let $T= \S^1 \times \R$ be the ``long tube'', and let $\mu: T \to T$ be the transformation defined by $\mu(z,x) = (-z,x+1)$.
We choose a small disk $D$ embedded in $\S^1 \times (0,1)$.
Then $\tilde N = T \setminus (\cup_{k \in \Z}\, \mu^k D)$ (see Figure 15\,(ii)), $N_{2,1}$ is the quotient $\tilde N/\langle \mu \rangle$, and $p$ is the quotient map $\tilde N \to \tilde N/\langle \mu \rangle = N_{2,1}$.

\bigskip\noindent
For $k \in \Z$ we set $\tilde a_k = \S^1 \times \{k\}$.
We can and do assume that $\tilde a_k$ is a lift of $a$, and that every lift of $a$ is of this form.
On the other hand, we denote by $\tilde d_k$ the boundary of $\mu^k\,D$.
Note that $\tilde d_k$ is a lift of $d$, and every lift of $d$ is of this form.

\bigskip\noindent
Let $b$ be a generic two-sided circle in $N_{2,1}$.
We can assume without loss of generality that $b$ intersects $a$ transversally in finitely many points. 
We say that $a$ and $b$ \emph{cobound a bigon} if there exists a disk $\Delta$ embedded in $N_{2,1}$ bounded by an arc of $a$ and an arc of $b$ (see Figure 15\,(iii)).
Note that, if such a bigon exists, then one can replace $b$ by a circle $b'$ isotopic to $b$ such that $|b' \cap a| = |b \cap a|-2$.

\bigskip\noindent
By Proposition 19, $b$ is isotopic to $a$ in $N_{2,0}$, hence all the lifts of $b$ in $\tilde N$ are circles.
Let $\tilde b_0$ be one of these lifts.
If $\tilde b_0 \cap \tilde a_k = \emptyset$ for all $k \in \Z$, then $\tilde b_0$ is isotopic to some $\tilde a_k$, hence $b$ is homotopic to $a$ (in $N_{2,1}$), and therefore, by Theorem 3, $b$ is isotopic to $a$.

\bigskip\noindent
So, we can suppose that there exists $k \in \Z$ such that $\tilde b_0 \cap \tilde a_k \neq \emptyset$.
Let $k_0,k_1 \in \Z$ such that $\tilde b_0 \cap \tilde a_{k_0} = \emptyset$, $\tilde b_0 \cap \tilde a_{k_0+1} \neq \emptyset$, $\tilde b_0 \cap \tilde a_{k_1} \neq \emptyset$, and $\tilde b_0 \cap \tilde a_{k_1+1} = \emptyset$ (see Figure 15\,(iv) and Figure~15\,(v)).
We have $k_0 < k_1$, $\tilde b_0 \cap \tilde a_\ell \neq \emptyset$ for all $\ell \in \{k_0+1, \dots, k_1\}$, and $\tilde b_0 \subset \S^1 \times (k_0,k_1+1)$.

\bigskip\noindent
There exists a bigon, $\Delta_0$, in $\S^1 \times (k_0,k_0+1]$, bounded by an arc of $\tilde b_0$ and an arc of $\tilde a_{k_0+1}$ (see Figure 15\,(iv)).
Similarly, there exists 	a bigon, $\Delta_1$, in $\S^1 \times [k_1,k_1+1)$, bounded by an arc of $\tilde b_0$ and an arc of $\tilde a_{k_1}$ (see Figure 15\,(v)).
If $\mu^{k_0}D$ is not contained in $\Delta_0$, that is, if $\Delta_0 \cap \mu^{k_0}D = \emptyset$, then $\Delta_0$ is contained in $\tilde N$ and $p(\Delta_0)$ is a bigon in $N_{2,1}$ bounded by an arc of $b$ and an arc of $a$.
As pointed out before, we can then reduce the number of intersection points between $a$ and $b$, and therefore apply an induction on $|a \cap b|$.

\bigskip\noindent
So, we can assume that $\mu^{k_0}D$ is contained in $\Delta_0$.
Similarly, we can assume that $\mu^{k_1}D$ is contained in $\Delta_1$.
Then $\mu^{k_1}D \subset \mu^{k_1-k_0} \Delta_0 \cap \Delta_1$, hence $\mu^{k_1-k_0} \Delta_0 \cap \Delta_1 \neq \emptyset$, and therefore $\mu^{k_1-k_0} \tilde b_0 \cap \tilde b_0 \neq \emptyset$ (since $\tilde a_{k_1} \cap \tilde a_{k_1+1} = \emptyset$): a contradiction.
\qed

\bigskip\noindent
Similar arguments apply to the solution of the following exercise. 

\bigskip\noindent
{\bf Exercise.}
Suppose that the torus, $S_{1,0}$, is obtained from the one-holed torus, $S_{1,1}$, by gluing a disk along the boundary component of $S_{1,1}$.
Prove that the map $\CC(S_{1,1}) \to \CC(S_{1,0})$, induced by the inclusion map $S_{1,1} \hookrightarrow S_{1,0}$, is a bijection. 

\bigskip\noindent
{\bf Definition.}
Assume that $N=N_{2,1}$ is the one-holed Klein bottle.
Consider the (oriented) circle $a$ and the crosscap $\mu$ depicted in Figure 16.
We orient $a$ and a regular neighborhood of $a$ as in the figure. 
We denote by $y_{\mu,a}$ the mapping class in $\MM(N,\partial N)$ obtained sliding $\mu$ once along $a$.
Then $y_{\mu,a}$ is called the \emph{crosscap slide} of $\mu$ along $a$.
Note that $y_{\mu,a}$ reverses the orientation on a regular neighborhood of $a$, $y_{\mu,a}^2=t_c$, and $y_{\mu,a} t_a y_{\mu,a}^{-1} = t_a^{-1}$, where $c$ is the boundary component of $N$.

\begin{figure}[tbh]
\bigskip\centerline{
\setlength{\unitlength}{0.4cm}
\begin{picture}(7,5)
\put(0,0){\includegraphics[width=2.4cm]{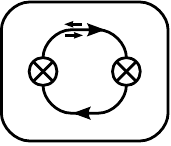}}
\put(0.5,3){\small $\mu$}
\put(4,4){\small $a$}
\put(6.2,4){\small $c$}
\end{picture}} 

\bigskip
\centerline{{\bf Figure 16.} Crosscap slide.}
\end{figure}

\bigskip\noindent
The proof of the following proposition is based on Proposition 21. 

\bigskip\noindent
{\bf Proposition 22}
(Stukow \cite{Stuko2}).
{\it $\MM(N_{2,1},\partial N_{2,1}) \simeq \Z \rtimes \Z$, where the first copy of $\Z$ is generated by $t_a$, and the second copy is generated by $y_{\mu,a}$.}

\bigskip\noindent
More generally, we have the following definition. 

\bigskip\noindent
{\bf Definition.}
Let $N$ be a non-orientable surface.
Let $\mu$ be a crosscap in $N$, and let $a$ be a generic two-sided (oriented) circle that intersects $\mu$ in a single point. 
Then the mapping class in $\MM(N)$ obtained by sliding $\mu$ once along $a$ is called the \emph{crosscap slide} of $\mu$ along $a$, and is denoted by $y_{\mu,a}$.

\bigskip\noindent
Now, we come back to our initial hypothesis.
So, $M$ denotes a surface which may be orientable or not, have boundary or not, and have more than one connected component.
However, we assume now that each connected component has a negative Euler characteristic. 
Recall that $\CC(M)$ denotes the set of isotopy classes of generic circles, and that $\CC_0(M)$ denotes the set of isotopy classes of generic two-sided circles. 

\bigskip\noindent
{\bf Definition.}
Recall that the \emph{intersection index} of two classes, $\alpha, \beta \in \CC(M)$, is $i(\alpha, \beta) = \min\{|a \cap b| \mid a \in \alpha \text{ and } b \in \beta\}$.
The set $\CC(M)$ is endowed with a structure of simplicial complex, where a simplex is a finite non-empty subset $\AA$ of $\CC(M)$ satisfying $i(\alpha, \beta)=0$ for all $\alpha, \beta \in \AA$.
This simplicial complex is called the \emph{curve complex} of $M$.
Note that a simplex $\AA$ admits a set of representatives, $A$, such that $a \cap b = \emptyset$ for all $a,b \in A$, $a \neq b$.
Such a set of representatives is called an \emph{admissible set of representatives} for $\AA$.
Note also that the group $\MM(M)$ acts naturally on $\CC(M)$ and on $\CC_0(M)$, and that these actions are simplicial. 

\bigskip\noindent
{\bf Definition.}
Let $\AA$ be a simplex in $\CC(M)$, and let $A$ be an admissible set of representatives for $\AA$.
We denote by $\MM_\AA(M)$ the stabilizer of $\AA$ in $\MM(M)$.
On the other hand, we denote by $M_\AA$ the natural compactification of $M \setminus (\cup_{a \in A} a)$.
The groups $\MM_\AA(M)$ and $\MM(M_\AA)$ are related by a homomorphism $\Lambda_\AA : \MM_\AA(M) \to \MM(M_\AA)$ defined as follows.
Let $f \in \MM_\AA(M)$.
Choose a representative $F \in \Homeo(M)$ of $f$ such that $F(A)=A$.
The restriction of $F$ to $M \setminus (\cup_{a \in A} a)$ extends in a unique way to a homeomorphism $\hat F \in \Homeo(M_\AA)$.
Then $\Lambda_\AA(f)$ is the element of $\MM(M_\AA)$ represented by $\hat F$.
The homomorphism $\Lambda_\AA$ is called the \emph{reduction homomorphism} along $\AA$.

\bigskip\noindent
{\bf Definition.}
For a simplex $\AA$ in $\CC(M)$, we set $\AA_0 = \AA \cap \CC_0(M)$, and we denote by $Z_\AA$ the subgroup of $\MM(M)$ generated by $\{ t_\alpha \mid \alpha \in \AA_0\}$.

\bigskip\noindent
The following proposition is a ``classical result'' for mapping class groups of connected orientable surfaces (see \cite{ParRol1}, for instance), and it is due to Stukow \cite{Stuko2, Stuko1} for the mapping class groups of connected non-orientable surfaces.
The general case is a straightforward consequence of these two cases. 

\bigskip\noindent
{\bf Proposition 23.}
{\it Let $\AA$ be a simplex of $\CC(M)$, and let $\Lambda_\AA : \MM_\AA (M) \to \MM(M_\AA)$ be the reduction homomorphism along $\AA$.
Then the kernel of $\Lambda_\AA$ is $Z_\AA$.}

\bigskip\noindent
{\bf Definition.}
Let $\Gamma(M)= \{ M_i \mid i \in I\}$ denote the set of connected components of $M$.
Let $f \in \MM(M)$.
We say that $f$ is \emph{pseudo-Anosov} if $\CC(M_i) \neq \emptyset$ for all $i \in I$, and  $f^n(\alpha) \neq \alpha$ for all $\alpha \in \CC(M)$ and all $n \in \Z \setminus \{ 0\}$.
We say that $f$ is \emph{reducible} if there exists a (non-empty) simplex $\AA$ of $\CC(M)$ such that $f(\AA)=\AA$.
In this case, the simplex $\AA$ is called a \emph{reduction simplex} for $f$, and an element of $\AA$ is called a \emph{reduction class} for $f$.

\bigskip\noindent
Here is an important property of pseudo-Anosov elements.
It is implicitly proved in \cite{FaLaPo1} for orientable surfaces, and it is a consequence of \cite{Wu1} (see Proposition 27 below) for non-orientable surfaces. 

\bigskip\noindent
{\bf Proposition 24.}
{\it Assume that $M$ is a connected surface. 
Let $f \in \MM(M)$ be a pseudo-Anosov element, and let $Z_f(\MM(M))$ be the centralizer of $f$ in $\MM(M)$.
Then $f$ has infinite order, and $\langle f \rangle$ has finite index in $Z_f(\MM(M))$.}

\bigskip\noindent
{\bf Definition.}
Recall that $\Gamma(M)= \{ M_i \mid i \in I\}$ denotes the set of connected components of $M$.
Let $\SSS(\Gamma(M))$ denote the group of permutations of $\Gamma(M)$.
Note that there is a natural homomorphism $\varphi : \MM(M) \to \SSS(\Gamma(M))$ whose kernel is naturally isomorphic to $\oplus_{i \in I} \MM(M_i)$.
Let $f \in \MM(M)$, and let $n \in \Z\setminus \{0\}$ such that $f^n \in \Ker\, \varphi$. 
We say that $f \in \MM(M)$ is \emph{adequately reduced} if, for all $i \in I$, the restriction of $f^n$ to $\MM(M_i)$ is either of finite order, or pseudo-Anosov. 
It is easily proved that this definition does not depend on the choice of $n$.
A reduction simplex $\AA$ for $f \in \MM(M)$ is called an \emph{adequate reduction simplex} for $f$ if the reduction of $f$ along $\AA$ is adequately reduced. 

\bigskip\noindent
The following result is extremely important in the study of mapping class groups of orientable surfaces. 

\bigskip\noindent
{\bf Theorem 25}
(Thurston \cite{FaLaPo1}).
{\it Assume that $M$ is an orientable surface. 
Then a mapping class $f \in \MM(M)$ is either reducible, or adequately reduced. 
Moreover, if $f$ is reducible, then there exists an adequate reduction system for $f$.}

\bigskip\noindent
The result of Theorem 25 holds for non-orientable surfaces, but I do not know whether the proof of \cite{FaLaPo1} can be extended to non-orientable surfaces.
However, following \cite{Wu1}, the result for non-orientable surfaces can be reduced to the result for orientable surfaces using oriented two-sheeted covers as follows. 

\bigskip\noindent
From now on, $N$ denotes a compact non-orientable surface possibly with boundary. 
We denote by $p : M \to N$ the two-sheeted cover of $N$ with $M$ orientable.
(Such a cover exists and is unique.) 
Let $J: M \to M$ be the (unique) non-trivial deck transformation of $p$.
We can and do assume that $N$ and $M$ are endowed with hyperbolic metrics such that the boundary components are geodesics, and $J$ is an isometry. 
We denote by $j$ the element of $\MM(M)$ represented by $J$.

\bigskip\noindent
{\bf Definition.}
Let $F\in \Homeo (N)$.
Then $F$ has two lifts on $M$, one which preserves the orientation, that we denote by $\tilde F$, and one which reverses the orientation, $J \circ \tilde F$.
These two lifts commute with $J$.
By \cite{Ziesc1} (see also \cite{BirChi1, BirHil1}), if $F$ and $F'$ are two homeomorphisms of $N$ such that $\tilde F$ and $\tilde F'$ are isotopic, then $F$ and $F'$ are also isotopic.
Let $\MM^+(M)$ denote the subgroup of $\MM(M)$ formed by the isotopy classes that preserve the orientation, and let $\ZZ_j^+(M)$ denote the subgroup of $\MM^+(M)$ formed by the mapping classes that commute with $j$.
Then, by the above, we have the following.

\bigskip\noindent
{\bf Proposition 26.}
{\it The map which sends $F \in \Homeo(N)$ to $\tilde F$ induces an injective homomorphism $\iota : \MM(N) \to \ZZ_j^+(M)$.}

\bigskip\noindent
One of the main results in \cite{Wu1} is the following.

\bigskip\noindent
{\bf Proposition 27}
(Wu \cite{Wu1}).
{\it Let $f \in \MM(N)$.
Then:
\begin{itemize}
\item[(1)]
$f$ is reducible if and only if $\iota(f)$ is reducible.
\item[(2)]
$f$ has finite order if and only if $\iota(f)$ has finite order.
\item[(3)]
$f$ is pseudo-Anosov if and only if $\iota(f)$ is pseudo-Anosov.
\end{itemize}}

\bigskip\noindent
Now, to get the full result of Theorem 25 for non-orientable surfaces, we need also to understand the reduction simplices. 

\bigskip\noindent
{\bf Definition.}
Let $\alpha \in \CC(N)$.
Choose a circle $a$ such that $\alpha = [a]$.
If $\alpha$ is one-sided, then $p^{-1}(a)$ is formed by a unique essential circle, $\tilde a$.
Then we set $p^{-1}(\alpha)= \{[\tilde a]\}$. 
If $\alpha$ is two-sided, then $p^{-1} (a)$ is formed by two disjoint essential circles, $\tilde a_1$ and $\tilde a_2$.
Then we set $p^{-1}(\alpha) = \{ [\tilde a_1], [\tilde a_2] \}$. 
It is easily checked that the definition of $p^{-1} (\alpha)$ does not depend on the choice of the representative $a$.
For a subset $\AA$ of $\CC(N)$, we set $p^{-1}(\AA) = \cup_{\alpha \in \AA} p^{-1}(\alpha)$.

\bigskip\noindent
{\bf Definition.}
Let $\tilde a$ be an essential circle in $M$. 
If either $J(\tilde a)=\tilde a$, or $\tilde a \cap J(\tilde a)=\emptyset$, then $p(\tilde a)$ is an essential circle in $N$.
Let $\tilde \alpha \in \CC(M)$.
Suppose that either $j(\tilde\alpha) = \tilde\alpha$, or $i(\tilde\alpha, j(\tilde\alpha))=0$.
Choose $\tilde a \in \tilde \alpha$ such that either $J(\tilde a) = \tilde a$, or $J(\tilde a) \cap \tilde a=\emptyset$ (take a geodesic which represents $\tilde \alpha$, for instance), and denote by $p(\tilde \alpha)$ the element of $\CC(N)$ represented by $p(\tilde a)$.
It is easily checked that $p(\tilde \alpha)$ is well-defined (namely, does not depend on the choice of $\tilde a$). 
Let $\tilde \AA$ be a simplex of $\CC(M)$ satisfying $j(\tilde \AA) = \tilde \AA$.
Then, for all $\tilde \alpha \in \tilde \AA$, we have either $j(\tilde \alpha) = \tilde \alpha$, or $i(\tilde \alpha, j(\tilde \alpha))=0$.
It follows that $p(\tilde \AA) = \{ p(\tilde \alpha) \mid \tilde \alpha \in \tilde \AA\}$ is a well-defined subset of $\CC(N)$.

\bigskip\noindent
The following lemma is easily proved choosing geodesics in the sets of admissible representatives of $\AA$ (in (1)), and of $\tilde \AA$ (in (2)).

\bigskip\noindent
{\bf Lemma 28.}
{\it \begin{itemize}
\item[(1)]
If $\AA$ is a simplex of $\CC(N)$, then $p^{-1} (\AA)$ is a simplex of $\CC(M)$.
\item[(2)]
If $\tilde \AA$ is a simplex of $\CC(M)$ such that $j(\tilde \AA) = \tilde \AA$, then $p(\tilde \AA)$ is a simplex of $\CC(N)$.
\end{itemize}}

\bigskip\noindent
Now, the rest of the statement of Theorem 25 for non-orientable surfaces follows from the next proposition. 

\bigskip\noindent
{\bf Proposition 29}
(Wu \cite{Wu1}).
{\it Let $f \in \MM(N)$, and let $\AA$ be a simplex of $\CC(N)$. 
\begin{itemize}
\item[(1)]
$\AA$ is a reduction simplex for $f$ if and only if $p^{-1} (\AA)$ is a reduction simplex for $\iota(f)$.
\item[(2)]
$\AA$ is an adequate reduction simplex for $f$ if and only if $p^{-1} (\AA)$ is an adequate reduction simplex for $\iota(f)$.
\end{itemize}}

\bigskip\noindent
{\bf Acknowledgement.}
I warmly thank Christian Bonatti who helped me to find the proof of Proposition 21. 

%%%%%%%%%%

%%%%%

\bigskip\bigskip\noindent
{\bf Luis Paris,}

\smallskip\noindent 
Université de Bourgogne, Institut de Mathématiques de Bourgogne, UMR 5584 du CNRS, B.P. 47870, 21078 Dijon cedex, France.

\smallskip\noindent
E-mail: {\tt lparis@u-bourgogne.fr}

%%%%%
\end{document}